\documentclass{tran-l}
\newtheorem{thm}{Theorem}
\newtheorem{cor}{Corollary}
\newtheorem{lem}{Lemma}
\newtheorem{prop}{Proposition}
\theoremstyle{definition}
\newtheorem{defn}{Definition}
\theoremstyle{remark}

\begin{document}
\title[The Embedding Theorem for a Mixed Verbal Wreath Product]{The Embedding
Theorem for a Mixed Verbal Wreath Product in Lie Algebras
Representations}
\author[L. A. Simonian]{L.~A.~Simonian}
\email{siml@rau.lv} \maketitle \vspace{-0.4in}
\begin{center}\small{
R$\bar{\i}$ga Aviation University, 1 Lomonosov St., R$\bar{\i}$ga
LV 1019, Latvia}\end{center}
\begin{abstract} The mixed verbal wreath product in
representations of Lie algebras defined in Reference~\cite{13} is
a construction parallel to the verbal wreath product of Lie
algebras introduced by A.~L.~Shmelkin~\cite{3}. It is shown that
the analog of the group theorem on embedding in the verbal wreath
product~\cite{3}, also holds for representations of Lie algebras.

{\sc Theorem}. {\sl Let $\mathcal X$ be a multihomogeneous variety
of representations of Lie algebras over a field $K$, and $L$ be an
absolutely free Lie algebra over $K$ with a system of free
generators$\{x_i,i\in I\}$. Suppose that $M$ is an ideal of $L$,
and $\mathcal F$ is absolutely free associative algebra with
unity, generated freely by $\{x_i,i\in I\}$. Furthermore, assume
that ${\mathcal F}_1$ is the universal enveloping algebra of $M$,
${\mathcal F}_1=U(M)$, $L_1$ is an absolutely free Lie algebra
over $K$ with a system of free generators $\{y_i,i\in I\}$, $\Phi$
is an absolutely free associative algebra with unity generated
freely by $\{y_i,i\in I\}$. Next, $(\Phi/{\mathcal
X}^*(\Phi,L_1),L_1)$ is assumed to be a free cyclic representation
of the variety $\mathcal X$ generated freely by $\{y_i,i\in I\}$.
Finally, let us denote by $\bar x$ the image of $x\in L$ in $L/M$
under canonical epimorphism. \emph{Then} the mapping $x_i\to
y_i+\bar x_i$ is extended to the monomorphism of $({\mathcal
F}/{\mathcal FX^*(F}_1,M),L)$ in $(\Phi/{\mathcal
X}^*(\Phi,L_1),L_1)\stackrel{\textstyle\mathcal X}{wr}L/M$.}

{\sc Corollary}. {\sl Let $\mathcal X$ be a multihomogeneous
variety of Lie algebras representation, and $\Theta$ be a variety
of Lie algebras over $K$. Then the free cyclic pair $({\mathcal
F/(X}\times\Theta)^*({\mathcal F},L),L)$ of the variety ${\mathcal
X}\times\Theta$ is isomorphic to subrepresentation of mixed wreath
$\mathcal X$-product $(\Phi/{\mathcal
X}^*(\Phi,L_1),L_1)\stackrel{\textstyle\mathcal X}{{wr}}L/
\Theta(L)$, generated by the elements $y_i+\bar x_i$. Here
$L,{\mathcal F},L_1$ and $\Phi$ are the same as in the Theorem and
$\bar x_i$ is the image of $x_i\in L$ in $L/\Theta(L)$.}
\end{abstract}
\vspace{0.5in}

The mixed verbal wreath product was introduced in \cite{13}. This
construction is parallel to verbal wreath product of Lie algebras
that had been defined by A. L. Shmelkin in \cite{3}. In this paper
we prove the embedding theorem for mixed verbal wreath products.
Let us give the definition of the mixed verbal wreath product (see
Ref.\cite{13}).
\begin{defn}
Let $\mathcal X$ be a variety of representations of Lie algebras
over a ring $K$ with 1, $(V,A)\in\mathcal X$, and $B$ be an
arbitrary Lie algebra over $K$. Lie pair $(G,\Gamma)$ is referred
to as {\it a mixed verbal wreath product} or {\it mixed wreath
$\mathcal X$-product} of $(V,A)\in\mathcal X$ and $B$, if the
following conditions are satisfied:
\begin{enumerate}
\item $(V,A)$ is a subpair of $(G,\Gamma)$;
\item $A$ and $B$ generate $\Gamma$;
\item $V$ generates $G$ as $\Gamma$-module;
\item If $P$ is ideal of $\Gamma$, generated by $A$, then pair
$(G,P)\in\mathcal X$;
\item If $\varphi:(V,A)\to(G_1,\Gamma_1)$ and $\psi:B\to\Gamma_1$
are homomorphisms such that
\begin{enumerate}
\item $\varphi(A)$ and $\psi(B)$ generate $\Gamma_1$,
\item $\varphi(V)$ generates $G_1$, as $\Gamma_1$-module.
\item $(G_1,\Phi)\in\mathcal X$, where $\Phi$ is ideal of $\Gamma_1$,
generated by $\varphi(A)$,
\end{enumerate}then the unique homomorphism
$\chi:(G,\Gamma)\to(G_1,\Gamma_1)$ exists, such that the restriction
of $\chi$ to $(V,A)$ coincides with $\varphi$,
and the restriction of $\chi$ to $B$ coincides with $\psi$.
\end{enumerate}
\end{defn}

We shall denote the mixed wreath $\mathcal X$-product of $(V,A)$
and $B$ by $(V,A)\stackrel{\textstyle\mathcal X}{wr}B$. Other
definitions of mixed verbal product and its the most important
properties are considered in \cite{13}. If $K$ is a field and
$(V,A)$ is a free pair of a multihomogeneous variety $\mathcal X$
then the mixed wreath $\mathcal X$-product may be constructed as
follows. Let  $Y=\{y_i,i\in I\}$ be a system of free generators of
$V$, $U=U(B)$ be the universal enveloping algebra for $B$,
$Z=\{z_n,n\in J\}$ be a basis of $U$ and $W=U(P)$ be the universal
enveloping algebra for $P$. Let us construct an $\mathcal X$-free
representation of Lie algebra $P$ in vector space generated by
elements $x_i\otimes z_n,i\in I,n\in J$, \emph{freely} with
respect to $P$. This representation will have the form:
\begin{equation}
((V^*\otimes_KU)\otimes_KW/{\mathcal
X}^*((V^*\otimes_KU)\otimes_KW),P),P) , \label{form}
\end{equation}
where $V^*$ is a subspace of $V$ spanned by $Y$, ${\mathcal
X}^*((V^*\otimes_KU)\otimes_KW),P)$ is a verbal submodule
corresponding to $\mathcal X$ of pair
$((V^*\otimes_KU)\otimes_KW),P)$. The action of $B$ on generators
$x_i\otimes z_n$ is determined by the rule: $(x_i\otimes
z_n)b=x_i\otimes z_nb$, where $z_nb$ is a product within $U$. Thus
$V^*\otimes_KU$ is free $U$-module generated freely by $y_i,i\in
I$.The action of $B$ on $V^*\otimes_KU$ may be unuquely extended
to the action on $(V^*\otimes_KU)\otimes_KW$ by the formula:
\begin{equation}
((v\otimes u)\otimes w)b=(v\otimes ub)\otimes w+(v\otimes u)\otimes[w,b].
\label{formula}
\end{equation}
Together with the action of $P$ on $(V^*\otimes_KU)\otimes_KW$ the
latter  determines the action of $\Gamma=P\lambda B$ on
$(V^*\otimes_K U)\otimes_KW$. Since the verbal submodule
\\${\mathcal X}^*((V^*\otimes_K U)\otimes_KW,P)$ is invariant
under derivations that are induced by derivations of $P$, we may
consider the following pair:

\begin{equation}
((V^*\otimes_KU)\otimes_KW/{\mathcal
X}^*((V^*\otimes_KU)\otimes_KW,P),\Gamma). \label{pair}
\end{equation}
This pair is the mixed wreath $\mathcal X$-product
$(V,A)\stackrel{\textstyle\mathcal X}{wr}B$.

A variety $\mathcal X$ is called a multihomogeneous variety if for
any identity $yv(x_1,
 x_2, \ldots,x_n) =0$ of variety $\mathcal X$
the identity $yv_{\alpha}(x_1,x_2, \ldots,x_n)
=0$ also holds in
$\mathcal X$. Here
$v=\sum\limits_{\alpha}v_{\alpha}(x_1,x_2,\ldots,x_n)$ is the
multihomogeneous decomposition of $v(x_1,x_2,\ldots,x_n)$.
\begin{prop}
 Let $(V,L)$ be a free cyclic representation of a
multihomogeneous variety $\mathcal X$ with the system of free
generators $X=\{x_\alpha,\alpha\in I\}$, $Y$ be a subset of $L$,
that is linearly independent modulo $L^2$ and $(W,M)$ be the
cyclic subrepresentation of representation $(V,L)$, generated by
$Y$ and $e$, where $e$ is the generator of $L$-module $V$. Then
$(W,M)$ is a free representation of the variety  $\mathcal X$.
\end{prop}
\begin{proof}
 Let $\mathcal S$ be a variety of all trivial representations
over field $K$, that are representations in which $L$ acts as zero
algebra. For $\mathcal X=S$ the assertion of the Proposition is
evident. Therefore, let us assume that $\mathcal X\ne S$. By the
condition $L$ is free Lie algebra with the system of free
generators $X=\{x_\alpha,\alpha\in I\}$, and $V={\mathcal F}/U$,
where $\mathcal F$ is free associative algebra with unity, that is
generated by $X$, and $U$ is the ideal of identities which respond
to $\mathcal X$. Here as a free generator of $V$ we take the image
of unity $e$ of the algebra $\mathcal F$ in $V$. Let $Z$ be a set
of the same cardinality as $Y$, and $\varphi$ is bijective map of
$Z$ onto $Y$. Let $(V_1,L_1)$ denote the free cyclic
representation of variety $\mathcal X$, generated by $Z$. The map
$\varphi$ can be extended to the homomorphism $\bar\varphi$ of
pair $(V_1,L_1)$ onto $(W,M)$. If $\bar\varphi$ is not
isomorphism, then for some nonzero element $v_1\in V_1$ we have
$\bar\varphi(v_1)=0$. Let $v_1$ be expressed in terms of the
elements $z_1,z_2,\ldots,z_m$ from $Z$, and let
$\varphi(z_i)=y_i,i=1,2,\ldots,m$. Denote by $(A,H)$ and $(B,G)$
the subpairs in $(V_1,L_1)$ and $(W,M)$ generated  by
$z_1,z_2,\ldots,z_m$ and $y_1,y_2,\ldots,y_m$, respectively. The
restriction $\bar\varphi$ to $(A,H)$ is homomorphism of the (free
in $\mathcal X$) pair $(A,H)$ onto $(B,G)$ and is not injective.
Therefore,  it is possible to limit ourselves to the case of
finite set $Y=\{y_1,y_2,\ldots,y_m\}$. Let
$\{x_1,x_2,\ldots,x_n\}$ be the set of elements from $X$, and let
the elements $\{y_1,y_2,\ldots,y_m\}$ be expressed through the
former. Denote by $L_1$ the subalgebra of $L$ generated by
$\{x_1,x_2,\ldots,x_n\}$. Since $L_1^2\subset L^2$, the elements
$y_1,y_2,\ldots,y_m$ are linearly independent modulo $L_1^2$.
Therefore it is possible to consider $X$ a finite set as well.
>From the definition of $L$ it then follows that
$\dim L/L^2=n\ge m$. Since $y_1,y_2,\ldots,y_m$ are linearly
independent modulo $L^2$, the system $y_1,y_2,\ldots,y_m$ can be
completed to the basis modulo $L^2$ with elements of
$\{x_1,x_2,\ldots,x_n\}$. Upon renumbering, if necessary, added
elements from $X$ we end up with a basis set
$Y'={y_1,y_2,\ldots,y_m,x_{m+1},x_{m+2},\ldots,x_n}$. Denote by
$M_1$ the subalgebra of $L$, generated by $Y'$ and by $(W_1,M_1)$
the corresponding subpair in $(V,L)$. The map $\psi:x_1\to
y_1,\ldots,x_m\to y_m,x_{m+1}\to x_{m+1},\ldots,x_n\to x_n$ is
extended to homomorphism $\psi$ of pair $(V,L)$ onto pair
$(W_1,M_1)$. Let us prove that $\psi$ is isomorphism. Suppose the
contrary: there is $v\ne0$ that belongs to $Ker\psi$. Since
$\bigcap\limits_{n=1}^\infty{\mathcal F'}^n=0$, an $n$ exists,
such that $v\notin{\mathcal F'}^n$. Since $\psi({\mathcal
F'}^n)\subseteq{\mathcal F'}^n$, mapping $\psi$ induces
homomorphism $\bar\psi$ of pair $({\mathcal F}/U+{\mathcal
F'}^n,L)\approx({\mathcal F}/U_2+U_3+ \cdots+U_{n-1}+{\mathcal
F'}^n,L)$ into itself. Here $U_i$ is homogeneous component of the
degree $i$ of the ideal of identities $U$ and $\mathcal F'$ is
ideal of $\mathcal F$, generated by $X$. Since $({\mathcal
F}/U+{\mathcal F'}^n,L)$ is finitely stable representation, the
set $Y'$ is a system of generators. Therefore, $\bar\psi$ is
surjective mapping of finitedimensional linear space ${\mathcal
F}/U+{\mathcal F'}^n$ onto itself. Consequently, $\bar\psi$ is
isomorphism. Hence, $v\in{\mathcal F'}^n$, which brings us to a
contradiction. Hence, $\psi$ is isomorphism and $(W_1,M_1)$ is
$\mathcal X$-free representation.
\end{proof}
\begin{lem}\cite{3} Let $A$ be an abelian ideal of a Lie algebra
$B$ over ring $K$, $G=B/A$ and $G$ --- free $K$-modules. Then $B$
is isomorphic to the subalgebra (that we denote also by $B$) of a
semidirect product $S=C\lambda G$, where $C$ is some abelian Lie
algebra and $C\cap B=A$.
\end{lem}
\begin{proof} As is well known, $A$ can be regarded as $G$-module
if the action of $G$ in $A$ is defined by the formula $ag=[a,b]$,
where $b$ is an arbitrary element of the coset $g$. Proceeding
from a decomposition $B=A\oplus H$, where $H$ is a compliment of
the subspace $A$ with respect to $B$, we can define a system of
representatives by means of a one-to-one linear mapping
$\sigma:G\to H$, such that $\overline{g\sigma}=g$. Here  the coset
containing $h$  is denoted by $\bar h$. The bilinear mapping
$f(g_1,g_2)=[g_1^{\sigma},g_2^{\sigma}]-[g_1,g_2]^{\sigma}$ of
$B\times B$ into $A$ produces the set of elements $f(g_i,g_j)$ in
$A$ that is usually called a factor set. $H$ is a subalgebra if
and only if
 $f(g_i,g_j)$
is zero mapping. For a given $f$, the Lie algebra $B$, which is an
extension of $A$ by $B/A$, is defined as a linear space
$A\times G$ with a commutator:
\begin{equation}[(a_1,g_1),(a_2,g_2)]=(a_1g_2-a_2g_1+f(g_1,g_2),[g_1,g_2]).
\label{commutator}
\end{equation}
The split extension is obtained when $f(g_1,g_2)\equiv 0$.

Let $U$ be a universal enveloping algebra of $G$; then $A$ is a
$U$-module.
$U$-module $A$ can be inserted isomorphically into some injective
$U$-module $C$. Let $S$ be an extension of $C$ by $G$ with the
factor set $f$. Such an extension is possible because $A$ is inserted
into $C$. It is clear that $B\subset S$ and $B\cap C=A$. Since $C$ is
the injective $G$-module, the extension $S$ is split\cite{18}.
Hence, $S$ can be considered as a semidirect product $C\lambda G$,
which makes the proof complete.
\end{proof}
Let us also define the isomorphism $\mu$ of the extension $S$
and the semidirect
product $C\lambda G$. As is well known (see Ref.\cite{82})
the extension $S$
is the split extension if and only if there is a linear mapping
$\rho$ of $G$ into $C$, such that
\begin{equation}
f(g_1,g_2)=g_1^\rho g_2-g_2^\rho g_1-[g_1,g_2]^\rho.
\label{split_extension}
\end{equation}
Let us define a mapping $\mu$ of $S$ onto $C\lambda G$ by
$(a,g)\mu=(a+g^\rho,g)$

\begin{lem} The mapping $\mu$ is an isomorphism.\end{lem}

\begin{proof} If $(a_1,g_1),(a_2,g_2)\in S$, then
$(a_1,g_1)\mu=(a_1+g_1^\rho,g_1),(f_2,g_2)\mu=(a_2+g_2^\rho,g_2)$,
and
$[(a_1,g_1)\mu,(a_2,g_2)\mu]=[(a_1+g_1^\rho,g_1),(a_2+g_2^\rho,g_2)]=
(a_1g_2-a_2g_1+g_1^\rho  g_2-g_2^\rho g_1,[g_1,g_2])$. On the
other hand, $[(a_1,g_1),(a_2,g_2)]\mu=
(a_1g_2-a_2g_1+f(g_1,g_2),[g_1,g_2])\mu=(a_1g_2-a_2g_1+
f(g_1,g_2)+[g_1,g_2]^\rho,[g_1,g_2])$. By the splitting condition
$(a_1g_2-a_2g_1+g_1^\rho  g_2-g_2^\rho g_1,[g_1,g_2])=
(a_1g_2-a_2g_1+f(g_1,g_2)+[g_1,g_2]^\rho,[g_1,g_2])$. Hence $\mu$
is the homomorphism $S$ into $C\lambda G$. Furthermore,
$(a,g)\mu=(a+g^\rho,g)=0$ implies $g=0$. Consequently, $g^\rho=0$
and $a=0$. Thus the mapping $\mu$ is the injection. The surjective
property of $\mu$ is no less obvious: for an arbitrary $(a,g)\in
C\lambda G$,  $(a-g^\rho,g)\mu=(a,g)$ is true. Thus $\mu$ is an
isomorphism.\end{proof}

In view of this, for the injection $\nu$ of $B$ into $S$
determined in the Lemma above, we have $b\nu=\bar b+c$,
where $b$ is any element of $B$ and $\bar b$ is its image
under epimorphism $B\to G$, and $c\in C$.

\begin{lem} Let $L=L(\{x_i,i\in I\})$ be an absolutely free Lie algebra
over ring $K$ and $\{x_i,i\in I\}$ be the set of its free
generators. Let $M$ be an ideal of $L$, $L_1=L(\{y_i,i\in I\})$ an
absolutely free Lie algebra with a set of free generators
$\{y_i,i\in I\}$, $L/M$ be a free $K$-module, and $P$ be an ideal
in $L_1*L/M$ generated by $L_1$. Then there is an injection of
$M/M^2$ into $P/P^2$.
\end{lem}
\begin{proof} $P$ is absolutely free algebra with a system of free
generators $\{y_ie_{\alpha_1}e_{\alpha_2}\cdots e_{\alpha_n},
i\in
I,\alpha_1\le\alpha_2\le\cdots\le\alpha_n\}$, where
$\{e_\alpha,\alpha\in A\}$ is a basis for $L/M$ over $K$, and $A$
is \emph{well-ordered set}\cite{3}. $L_1*L/M$ is semidirect
product $P\lambda L/M$. The mapping $\varphi:x_i\to y_i+\bar x_i$
-- where $\bar x_i$ is the image of $x_i$ in $L/M$ -- is extended
to homomorphism $\varphi:L\to L_1*L/M$, and $\varphi$ maps $M$
into $P$. Since $P^2$ is an ideal of $P\lambda L/M$, it is
possible to consider the algebra $(P/P^2)\lambda L/M$ and
canonical homomorphism $\kappa:L_1*L/M\to(P/P^2)\lambda L/M$. The
composition $\chi$ of homomorphisms $\varphi$ and $\kappa$ is the
homomorphism of $L$ into $(P/P^2)\lambda L/M$. The kernel of
$\chi$ contains $M^2$ and, on the other hand, lies in $M$. The
homomorphism $\chi$ induces homomorphism
$\tilde\chi:L/M^2\to(P/P^2)\lambda L/M$. The algebra $L/M^2$,
being the extension of $M/M^2$ by $L/M$, can be isomorphically
embedded, by Lemma 2, into some semidirect product $S=C\lambda
L/M$. Let us denote this injection by $\nu$. Then by lemma 1, we
have $\nu(x_i)=c_i+\bar x_i$. Here again we have denoted by $x_i$
the image of $x_i$ in $L/M^2$ and $c_i$ is the same element of
$C$. By definition of the free product there is a homomorphism
$\sigma:L_1*L/M\to C\lambda L/M$ extending the mapping $\gamma$ of
$L_1$ into $C$ (that is defined by the formula $y_i\gamma=c_i$)
and the identity mapping of $L/M$ onto itself. Since $C$ is
abelian the ideal $P^2$ lies in the kernel of this mapping. Thus
we have the mapping $\sigma:(P/P^2)\lambda L/M$ into $C\lambda
L/M$ and, at the same time, the commutative diagram
$$\begin{array}{ccccc}L/M^2&\stackrel{\tilde\chi}{\longrightarrow}&
(P/P^2)\lambda L/M\\{}\\ \downarrow\lefteqn{\nu}&\swarrow\lefteqn
{\tilde\sigma}\\{}\\C\lambda L/M\end{array}.$$ Since $\nu$ is a
monomorphism and $\nu=\tilde\sigma\tilde\chi$, the mapping
$\tilde\chi$ is also monomorphism and, consequently, $M/M^2$ is
isomorphically embedded into $P/P^2$.
\end{proof}
It should be noted that Lemma~3 is essentially a parafrase of the
corresponding assertion in Ref.\cite{3}. For the discussion that
follows, we need this result in the specific form of Lemma~3.

\begin{thm} Let $\mathcal X$ be a multihomogeneous variety of representations
of Lie algebras over a field $K$, and $L$ be an absolutely free
Lie algebra over $K$ with a system of free generators$\{x_i,i\in
I\}$. Suppose that $M$ is an ideal of $L$, and $\mathcal F$ is
absolutely free associative algebra with unity, generated freely
by $\{x_i,i\in I\}$. Furthemore, assume that ${\mathcal F}_1$ is
the universal enveloping algebra of $M$, ${\mathcal F}_1=U(M)$,
$L_1$ is an absolutely free Lie algebra over $K$ with a system of
free generators $\{y_i,i\in I\}$, $\Phi$ is an absolutely free
associative algebra with unity generated freely by $\{y_i,i\in
I\}$. Next, $(\Phi/{\mathcal X}^*(\Phi,L_1),L_1)$ is assumed to be
a free cyclic representation of the variety $\mathcal X$ generated
freely by $\{y_i,i\in I\}$. Finally, let us denote by $\bar x$ the
image of $x\in L$ in $L/M$ under canonical epimorphism.
\emph{Then} the mapping $x_i\to y_i+\bar x_i$ is extended to the
monomorphism of $({\mathcal F}/{\mathcal FX^*(F}_1,M),L)$ in
$(\Phi/{\mathcal X}^*(\Phi,L_1),L_1)\stackrel{\textstyle\mathcal
X}{wr}L/M$.\end{thm}

\begin{proof} The mapping $\varphi:x_i\to y_i+\bar x_i$ is extended to the
homomorphism $\varphi:L\to L_1*L/M$, with $\varphi$ mapping $M$
into $P$. The homomorphism $\varphi$ is extended to the
homomorphism $\varphi:{\mathcal F}\to U(L_1*L/M)$, with $\varphi$
mapping ${\mathcal F}_1$ into $W=U(P)$. The composition of
$\varphi$ and the canonical homomorphism of $U(L_1*L/M)$ onto
$U(L_1*L/M)/{\mathcal X}^*(U(L_1*L/M),P)$ yields a homomorphism
$\psi:{\mathcal F}\to U(L_1*L/M)/{\mathcal X}^*(U(L_1*L/M),P)$.
The homomorphism $\varphi$ maps ${\mathcal X^*(F}_1,M)$ into
${\mathcal X}^*(W,P)$, and therefore, ${\mathcal FX^*(F}_1,M)$ is
mapped into $U(L_1*L/M){\mathcal X}^*(W,P)={\mathcal
X}^*(U(L_1*L/M),P)$. Thus we have the homomorphism
$\tilde\psi:{\mathcal F/FX^*(F}_1,M)\to U(L_1*L/M)/{\mathcal
X}^*(U(L_1*L/M),P)$ that is simultaneously the homomorphism of the
representation $({\mathcal F/FX^*(F}_1,M),L)$ in
$(U(L_1*L/M)/{\mathcal X}^*(U(L_1*L/M),P),L_1*L/M)$. Now let's
make use of the fact that $L_1*L/M=P\lambda L/M$ and
$U(L_1*L/M)=U\otimes_KW$, where $U=U(L/M)$. But then
$(U(L_1*L/M)/{\mathcal X}^*(U(L_1*L/M),P),L*L/M)$ and
$(U\otimes_KW/{\mathcal X}^*(W,P),P\lambda L/M)$ are isomorphic.
We then label this isomorphism $\varepsilon_2$.

Consider now the representation $({\mathcal F/FX^*(F}_1,M),L)$.
$e_1,e_2,\ldots,e_n,\ldots$ is a well-ordered basis of a
complement of $M$ with respect to $L$. Then ${\mathcal F=F}_1+\sum
e_{i_1}e_{i_2}\cdots e_{i_n}{\mathcal F}_1, i_1\le i_2\le\cdots\le
i_n$. This implies ${\mathcal FX^*(F}_1,M)={\mathcal
X^*(F}_1,M)+\sum e_{i_1}e_{i_2} \cdots e_{i_n}{\mathcal
X^*(F}_1,M)$ and ${\mathcal F}_1\cap{\mathcal FX^*(F}_1,M)=
{\mathcal X^*(F}_1,M), e_{i_1}e_{i_2}\cdots e_{i_n}{\mathcal
F}_1\cap{\mathcal FX^*( F}_1,M)= e_{i_1}e_{i_2}\cdots
e_{i_n}{\mathcal X^*(F}_1,M)$.

Then ${\mathcal F/FX^*(F}_1,M)=({\mathcal F}_1+{\mathcal
FX^*(F}_1,M))/{\mathcal FX^*(F}_1,M)+\sum(e_{i_1}e_{i_2}\cdots
e_{i_n}{\mathcal F}_1+ {\mathcal FX^*(F}_1,M))/{\mathcal
FX^*(F}_1,M)\approx{\mathcal F}_1/{\mathcal X^*F}_1,M)+ \sum
e_{i_1}e_{i_2}\cdots e_{i_n}({\mathcal F}_1/{\mathcal
X^*(F}_1,M))$. Let us label this isomorphism by $\varepsilon_1$.
The isomorphism $\varepsilon_1$ is also the isomorphism of
representations $({\mathcal F/FX^*(F}_1,M),L)$ and $({\mathcal
F}_1/{\mathcal X^*(F}_1,M)
+ \sum e_{i_1}e_{i_2}
\cdots
e_{i_n}({\mathcal F}_1/{\mathcal X^*(F}_1,M)),L)$. The action of
$L$ in the latter pair is defined by the rule: if $f\in{\mathcal
F}_1$ and $\bar f$ is image of $f$ in the canonical homomorphism
${\mathcal F}_1\to{\mathcal F}_1/{\mathcal X^*(F}_1,M)$, then
$(e_{i_1}e_{i_2}\cdots e_{i_n}\bar f)m=e_{i_1}e_{i_2}\cdots
e_{i_n} \overline{fm}$, if $m\in M$ and $(e_{i_1}e_{i_2}\cdots
e_{i_n}\bar f)p= (e_{i_1}e_{i_2}\cdots e_{i_n}p)\bar
f+e_{i_1}e_{i_2}\cdots e_{i_n} \overline{[f,p]}$, if $p$ belongs
to  a complement of $M$ with respect to $L$, for which the basis
$e_1,e_2,\ldots e_n,\ldots$ has been chosen. Recall that the
action in a pair $(U\otimes_KW/{\mathcal X}^*(W,P),P\lambda L/M)$
is defined in a similar way, namely: if $w\in W$ and $\bar w$ is
an image of $w$ under the homomorphism $W\to W/{\mathcal
X}^*(W,P)$, then $(u\otimes\bar w)p=u\otimes\overline{wp}$, when
$p\in P$ and $(u\otimes\bar w)l=ul\otimes\overline{[w,l]}$, if
$l\in L/M$.

Let us consider the following commutative diagram
$$\begin{array}{ccc} ({\mathcal F/FX^*(F}_1,M),L) &
\stackrel{\tilde\psi}{\longrightarrow} & (U({\mathcal
L})/{\mathcal X}^*(U({\mathcal L}),P),{\mathcal L})\\
\downarrow\lefteqn{\varepsilon_1} &&
\downarrow\lefteqn{\varepsilon_2}\\ E &
\stackrel{\delta}{\longrightarrow} & (U\otimes_KW/{\mathcal
X}^*(W,P), P\lambda L/M)
\end{array},$$
where ${\mathcal L}=L_1*L/M$ and $E=({\mathcal F}_1/{\mathcal
X^*(F}_1,M)+ \sum e_{i_1}\cdots e_{i_n}{\mathcal F}_1/{\mathcal
X^*(F}_1,M),L)$, and prove that $\delta$ is an isomorphism. Since
$\varepsilon_1$ is a monomorphism the equality
$\delta\varepsilon_1=\tilde\psi\varepsilon_2$ implies that
$\tilde\psi$ is also a monomorphism. The mapping $\varphi$ maps
$M$ into $P$ and ${\mathcal F}_1$ into $W$. The composition of
$\varphi$ and canonical homomorphism $W\to W/{\mathcal X}^*(W,P)$
gives the homomorphism $\gamma:{\mathcal F}_1\to W/{\mathcal
X}^*(W,P)$ which, in turn, induces the homomorphism
$\tilde\gamma:{\mathcal F}_1/{\mathcal X^*(F}_1,P) \to W/{\mathcal
X}^*(W,P)$. By the definition of $\delta$, its restriction on
${\mathcal F}_1/{\mathcal X^*(F}_1,M)$ coincides with
$\tilde\gamma$. Indeed, if $f_1\in{\mathcal F}_1$, then
$\delta(f_1+{\mathcal X^*(F}_1,M))=\varepsilon_2\tilde\psi
\varepsilon_1^{-1}(f_1+ {\mathcal
X^*(F}_1,M))=\varepsilon_2\tilde\psi(f_1+{\mathcal FX^*(F}_1,M))=
\varepsilon_2(\varphi(f_1)+{\mathcal X}^*(U\otimes_K{\mathcal
X}^*(W,P)))= \varphi(f_1)+{\mathcal
X}^*(W,P)=\tilde\gamma(f_1+{\mathcal X^*(F}_1,M))$. Let us prove
that $\tilde\gamma$ is a monomorphism. By Lemma~3, the mapping
$\varphi(x_i)=y_i+\bar x_i$ induces injection of $M/M^2$ into
$P/P^2$. Therefore, by Proposition, the pair
$(\tilde\gamma({\mathcal F}_1/{\mathcal X^*(F}_1,M),\varphi(M))$
is $\mathcal X$-free pair, generated freely by elements of
$\varphi(M)$ that are linearly independent modulo $P^2$. This
implies that $\tilde\gamma$ is a monomorphism. The mapping
$\varphi$ maps $e_i$ to $\varphi(e_i)=w_i+\bar e_i$, where $w_i\in
P$, and $\bar e_i\in L/M$. It is clear that $\bar e_i$ forms a
basis of $L/M$. Now let us order the monomials
$e_{i_1}e_{i_2}\cdots e_{i_n}$ postulating that a monomial of a
greater degree is greater; the monomials of the same degree are
ordered lexicographically. This ordering rule is carried over to
monomials $\bar e_{i_1}\bar e_{i_2}\cdots\bar e_{i_n}$ that form a
basis of $U=U(L/M)$.

Now let $u=\bar f_1+\sum e_{i_1}e_{i_2}\cdots e_{i_n} \bar
f_{i_1i_2\cdots i_n}$ be an arbitrary element ${\mathcal
F}_1/{\mathcal X^*(F}_1,M)+\sum e_{i_1}e_{i_2} \cdots
e_{i_n}({\mathcal F}_1/{\mathcal X^*(F}_1,M))$. Here $\bar
f_1=f_1+{\mathcal X^*(F}_1,M), \bar f_{i_1i_2\cdots
i_n}=f_{i_1i_2\cdots i_n}+{\mathcal X^*(F}_1,M)$, and
$f_1,f_{i_1i_2\cdots i_n}\in{\mathcal F}_1$. We need to prove that
$\delta(u)\ne 0$. One has: $\delta(u)=\delta(\bar f_1+ \sum
e_{i_1}e_{i_2}\cdots e_{i_n}\bar f_{i_1i_2\cdots i_n})=
\varepsilon_2\tilde\psi\varepsilon_1^{-1}((f_1+{\mathcal
X^*(F}_1,M))+ \sum e_{i_1}e_{i_2}\cdots e_{i_n}(f_{i_1i_2\cdots
i_n}+{\mathcal X^*(F}_1,M)))=
\varepsilon_2\tilde\psi((f_1+{\mathcal FX^*(F}_1,M))+ \sum
(e_{i_1}e_{i_2}\cdots e_{i_n}f_{i_1i_2\cdots i_n}+{\mathcal
FX^*(F}_1,M)))= \varepsilon_2((\varphi(f_1)+{\mathcal
X}^*(U\otimes_K{\mathcal X}^*(W,P)))+ \sum(w_{i_1}+ \bar
e_{i_1})(w_{i_2}+\bar e_{i_2})\cdots(w_{i_n}+\bar e_{i_n})
\varphi(f_{i_1i_2\cdots i_n})+ {\mathcal X}^*(U\otimes_K{\mathcal
X}^*(W,P))))$. If we decompose $(w_{i_1}+\bar
e_{i_1})(w_{i_2}+\bar e_{i_2}) \cdots(w_{i_n}+ \bar
e_{i_n})\varphi(f_{i_1i_2\cdots i_n})$ by the components of the
sum $U\otimes_KW=W+\sum\bar e_{i_1}\bar e_{i_2}\cdots\bar
e_{i_n}W$, then $\bar e_{j_1}\bar e_{j_2}\cdots\bar
e_{j_s}\varphi(f_{j_1j_2 \cdots j_s})$ will be the leading
elements of this decomposition. Here $e_{j_1}e_{j_2} \cdots
e_{j_s}f_{j_1j_2\cdots j_s}$ is the leading element of $\sum
e_{i_1}e_{i_2}\cdots e_{i_n}f_{i_1i_2\cdots i_n}$. Next, we have
$\varepsilon_2(\bar e_{j_1}\bar e_{j_2}\cdots\bar e_{j_s}
\varphi(f_{j_1j_2\cdots j_s})+{\mathcal X}^*(U\otimes_K{\mathcal
X}^*(W,P))))= \bar e_{j_1}\bar e_{j_2}\cdots \bar
e_{j_s}(\varphi(f_{j_1j_2 \cdots j_s})+{\mathcal X}^*(W,P))= \bar
e_{j_1}\bar e_{j_2}\cdots\bar e_{j_s}\tilde\gamma(f_{j_1j_2 \cdots
j_s}+{\mathcal X^*(F}_1,M))$. Since $\tilde\gamma$ is
amonomorphism, $\tilde\gamma(f_{j_1j_2 \cdots j_s}+{\mathcal
X^*(F}_1,M))\ne0$. This implies $\bar e_{j_1}\bar
e_{j_2}\cdots\bar e_{j_s} \bar\gamma(f_{j_1j_2\cdots
j_s}+{\mathcal X^*(F}_1,M))\ne0$ and, consequently,
$\delta(u)\ne0$.
\end{proof}

\begin{cor} Let $\mathcal X$ be a multihomogeneous variety of Lie algebras
representation, and $\Theta$ be a variety of Lie algebras over
$K$. Then the free cyclic pair $({\mathcal
F/(X}\times\Theta)^*({\mathcal F},L),L)$ of the variety ${\mathcal
X}\times\Theta$ is isomorphic to subrepresentation of mixed
$\mathcal X$-wreath product $(\Phi/{\mathcal
X}^*(\Phi,L_1),L_1)\stackrel{\textstyle\mathcal X}{{wr}}L/
\Theta(L)$, generated by the elements $y_i+\bar x_i$. Here
$L,{\mathcal F},L_1$ and $\Phi$ are the same as in the Theorem and
$\bar x_i$ is the image of $x_i\in L$ in $L/\Theta(L)$.\end{cor}

\begin{proof} In the Theorem take $M=\Theta(L)$ and use $({\mathcal
X}\times\Theta)^*({\mathcal F},L)={\mathcal FX^*(F}_1,\Theta(L))$.
Here ${\mathcal F}_1=U(\Theta(L))$.\end{proof}

\begin{cor} If $I$ is an
infinite set then ${\mathcal X}\times\Theta=
var(\Phi/{\mathcal
X}^*(\Phi,L_1),L_1)\stackrel{\textstyle\mathcal X}
{{wr}}L/\Theta(L))$\end{cor}

\begin{proof}$(\Phi/{\mathcal X}^*(\Phi,L_1),L_1)\stackrel{\textstyle\mathcal X}
{wr}L/\Theta(L)$ belongs to ${\mathcal X}\times\Theta$. Since --
by the Corollary~1 -- $(\Phi/{\mathcal
X}^*(\Phi,L_1),L_1)\stackrel{\textstyle\mathcal X}{wr}L/
\Theta(L)$ contains a free pair ${\mathcal X}\times\Theta$, then
$(\Phi/{\mathcal X}^*(\Phi,L_1),L_1)\stackrel{\textstyle\mathcal
X}{wr}L/ \Theta(L)$ generates ${\mathcal
X}\times\Theta$.\end{proof}

\end{document}